\documentclass[12pt,leqno]{article}   %(or equivalent)
\usepackage{amsfonts}
\usepackage{amssymb}
\usepackage{enumerate}
\usepackage{amsmath}
\usepackage[T1]{fontenc}
\usepackage{graphicx}

\newtheorem{ttt}{Theorem}[section]
\newtheorem{llll}[ttt]{Lemma}
\newtheorem{ccc}[ttt]{Claim}
\newtheorem{eee}[ttt]{Example}
\newtheorem{fff}[ttt]{Fact}
\newtheorem{rrr}[ttt]{Remark}
\newtheorem{sss}[ttt]{Statement}
\newtheorem{ddd}[ttt]{Definition}
\newtheorem{qqq}[ttt]{Question}
\newtheorem{cccc}[ttt]{Corollary}
\newtheorem{nnn}[ttt]{Notation}

\newcommand{\bt}{\begin{ttt}}
\newcommand{\bl}{\begin{llll}}
\newcommand{\bc}{\begin{ccc}}
\newcommand{\bex}{\begin{eee}}
\newcommand{\bfa}{\begin{fff}}
\newcommand{\br}{\begin{rrr}\upshape}
\newcommand{\bst}{\begin{sss}}
\newcommand{\bd}{\begin{ddd}\upshape}
\newcommand{\bq}{\begin{qqq}}
\newcommand{\bnn}{\begin{nnn}\upshape}
\newcommand{\bcor}{\begin{cccc}}

\newcommand{\bp}{\noindent\textbf{Proof. }}

\newcommand{\et}{\end{ttt}}
\newcommand{\el}{\end{llll}}
\newcommand{\ec}{\end{ccc}}
\newcommand{\eex}{\end{eee}}
\newcommand{\efa}{\end{fff}}
\newcommand{\er}{\end{rrr}}
\newcommand{\est}{\end{sss}}
\newcommand{\ed}{\end{ddd}}
\newcommand{\eq}{\end{qqq}}
\newcommand{\ecor}{\end{cccc}}
\newcommand{\enn}{\end{nnn}}

\newcommand{\ep}{\hspace{\stretch{1}}$\square$\medskip}

\newcommand{\lab}[1]{\label{#1}}

\newcommand{\NN}{\mathbb{N}}

\newcommand{\RR}{\mathbb{R}}

\newcommand{\al}{\alpha}
\newcommand{\ga}{\gamma} 
\newcommand{\de}{\delta}

\newcommand{\om}{\omega}
\newcommand{\si}{\sigma}
\newcommand{\ka}{\kappa}
\newcommand{\la}{\lambda}
\newcommand{\las}{\lambda^*} %Lebesgue outer measure

\newcommand{\iB}{\mathcal{B}}
\newcommand{\iC}{\mathcal{C}}
\newcommand{\iD}{\mathcal{D}}
\newcommand{\iG}{\mathcal{G}}
\newcommand{\iH}{\mathcal{H}}
\newcommand{\iI}{\mathcal{I}}
\newcommand{\iL}{\mathcal{L}}

\newcommand{\iF}{\mathcal{F}}
\newcommand{\iP}{\mathcal{P}}

\newcommand{\iN}{\mathcal{N}}

\newcommand{\sm}{\setminus}
\newcommand{\eset}{\emptyset}

\newcommand{\beq}{\begin{equation}}
\newcommand{\eeq}{\end{equation}}
\def\su{\subseteq}

\def\Leb{\iL} %Lebesgue measurable sets

\newcommand{\CH}{\textit{CH}}

\newcommand{\ZFC}{\textit{ZFC}}

\newcommand{\non}{{\rm non}}
\newcommand{\cof}{{\rm cof}}
\newcommand{\add}{{\rm add}}

\title{Can we assign the Borel hulls in a monotone way?}

\author{M\'arton Elekes\thanks{Partially supported by Hungarian Scientific
Foundation grants no.~49786, 61600 and F~43620.} \ 
and 
Andr\'as M\'ath\'e\thanks{Partially supported by Hungarian Scientific
Foundation grants no.~T~49786 and T~72655.}}

\begin{document}

\maketitle 

\begin{abstract}
A \emph{hull} of $A\su[0,1]$ is a set $H$ containing $A$ such that
%$\la(H\cap I)=\la(A\cap I)$ for every Lebesgue measurable set $I$.
%the inner Lebesgue measure of $H\setminus A$ is zero.
$\las(H)=\las(A)$.
We investigate all four versions of the following problem. Does there exist a
monotone (wrt.~inclusion) map that assigns a Borel/$G_\de$ hull to every
negligible/measurable subset of $[0,1]$?

Three versions turn out to be independent of \ZFC\ (the usual Zermelo-Fraenkel
axioms with the Axiom of Choice), while in the fourth case we only prove that
the nonexistence of a monotone $G_\de$ hull operation for all
measurable sets is consistent. It remains open whether existence here is also
consistent. We also answer a question of Z.~Gyenes and D.~P\'alv\"olgyi which
asks if monotone hulls can be defined for every chain (wrt.~inclusion) of
measurable sets. We also comment on the problem of hulls of all subsets of
$[0,1]$.
\end{abstract}

\insert\footins{\footnotesize{MSC codes: Primary 28A51;
Secondary 03E15, 03E17, 03E35, 28A05, 28E15, 54H05. }} 
\insert\footins{\footnotesize{Key Words: hull, envelope, Borel, 
    %delta,
    monotone, Lebesgue, measure, Cohen real, Continuum Hypothesis, \CH, \add, \cof,
    \non, descriptive set theory}}

\section{Introduction}

Let us fix some notation before formulating the problems of this note. 

\bnn
Let us denote by $\iN, \Leb, \iB$ and $\iG_\de$ the class of
Lebesgue negligible, Lebesgue measurable, Borel and $G_\de$ subsets of $[0,1]$,
respectively. Let $\la$ stand for Lebesgue measure, and $\las$ for Lebesgue outer measure.
\enn

\bd
A set $H\su[0,1]$ is a \emph{hull} of $A\su [0,1]$, if
$A$ is a subset of $H$ and $\las(H)=\las(A)$.
%$A\su H$
%and $\la(H\cap I)=\la(A\cap I)$ for
%every Lebesgue measurable set $I$.
%and the inner Lebesgue measure of $H\setminus A$ is zero; that is, $\la(M)=0$ for every
%measurable set $M\subset H\setminus A$.
\ed

%Clearly, if $H$ is a hull of $A$, then $\las(A)=\las(H)$. Moreover, if $\las(A)<\infty$, then $H$ is a hull of $A$ if and only if $A\su H$ and $\las(A)=\las(H)$.

Clearly, every set has a Borel, even a $G_\de$ hull. It is then very natural to ask
whether `a bigger set has a bigger hull'. (For the two actual motivations of
this paper see below.)

\bd
Let $\iD$ and $\iH$ be two subclasses of $\iP([0,1])$ (usually $\iD$ is $\iN$
or $\Leb$, and $\iH$ is $\iB$ or $\iG_\de$). 
If there exists a map $\varphi : \iD \to \iH$ such that
\begin{enumerate}
\item $\varphi(D)$ is a hull of $D$ for every $D\in\iD$
\item $D\su D'$ implies $\varphi(D)\su\varphi(D')$ 
\end{enumerate}
then we say that \emph{a monotone $\iH$ hull operation on $\iD$ exists}. 
\ed

The four questions we address in this paper are the following.

\bq
\lab{q:main}
Let $\iD$ be either $\iN$ or $\Leb$, and let $\iH$ be either $\iB$ or
$\iG_\de$. Does there exist a monotone $\iH$ hull operation on $\iD$?
\eq

\br
\lab{r:rem}
\mbox{}\par
\begin{enumerate} 

\item
The problem was originally motivated by the following question of
Z.~Gyenes and D.~P\'alv\"olgyi \cite{GP}. 
Suppose that $\iC\su\Leb$ is a chain of sets,
i.e.~for every $C,C'\in\iC$ either $C\su C'$ or $C'\su C$ holds. 
Does there exist a monotone $\iB/G_\de$ hull operation on $\iC$?
%Can we assign
%$G_\de$ hulls in a monotone way to elements of $\iC$?

\item 
Another motivation for our set of problems is that it seems to be 
very closely related to the theory of so called \emph{liftings}.
A map $l:\Leb\to\Leb$ is called a lifting if it preserves $\eset$, finite unions
and complement, moreover, it is constant on the equivalence classes modulo
nullsets, and also it maps each equivalence class to one of its members. 
Note that liftings are clearly monotone. 
For a survey of this theory see the chapter by Strauss, Macheras and Musia{\l}
in \cite{HM}, or the chapter by Fremlin in \cite{HB}, or Fremlin \cite{Fr}.
Note that the existence of Borel liftings is known to be independent of
\ZFC, but the existence of a lifting with range in a fixed Borel class is not
known to be consistent. 

We also remark that liftings are usually considered as $l^*:\Leb/\iN\to\Leb$ or
$l^*:\iP([0,1])/\iN\to\Leb$ maps. 

\item
\lab{s}
In light of the theory of liftings it is natural to ask if a monotone
Borel/$G_\de$ hull operation on $\iP([0,1])$ (i.e. all subsets of $[0,1]$) can be
defined.
We will see in Section \ref{s:all} that this is actually equivalent to the
existence of a monotone Borel/$G_\de$ hull operation on $\Leb$.

\end{enumerate}
\er

\br\lab{r:Polish}
We can extend the notion of \emph{hull} to any 
uncountable Polish space endowed with a nonzero
continuous $\sigma$-finite Borel measure $\mu$. Let $\mu^*$ denote the corresponding outer measure. If $\mu$ is finite, then we can define $H$ to be a hull of $A$ if $$H\supseteq A \quad \text{and} \quad \mu^*(H)=\mu^*(A).$$ However, if $\mu$ is infinite, then we say that a set $H$ is a hull of $A$ if
$$H\supseteq A \quad \text{and} \quad \mu^*(H\cap I)=\mu^*(A\cap I)$$
for every $\mu$-measurable set $I$. This latter property is in fact equivalent to that $\mu(M)=0$ for every $\mu$-measurable set $M\su H\setminus A$.

%We remark here that throughout this paper $[0,1]$ could be replaced by $\RR$, or $\RR^n$,
%or more generally, 
We remark here that the results (and proofs) of this paper remain valid if we replace $[0,1]$ by $\RR$, or by $\RR^n$, or more generally,  by any
uncountable Polish space endowed with a nonzero continuous $\sigma$-finite Borel measure.
Statement~\ref{s:equiv} is still true in this more general setting, as one can combine Lemma~\ref{l:dens} with the fact the every such Polish space is Borel isomorphic
(with a measure preserving isomorphism) either to the real line, or to a subinterval $[a,b]$ of the real line \cite{Ke}.

%(Statement~\ref{s:equiv} is still true in this more general setting, as every such space is Borel isomorphic
%(with a measure preserving isomorphism) either to the real line, or to a subinterval $[a,b]$ of the real line.)
%%(The arguments using the density
%%topology can be got around using that for such measures there exists a measure
%%preserving Borel isomorphism with a subinterval of $\RR$ \cite{Ke}.)
\er

The paper is organized as follows. First, in the next section we settle the
independence of the existence of a monotone Borel/$G_\de$ hull on $\iN$. The
consistency of the nonexistence immediately yields the consistency of the
nonexistence of a monotone Borel/$G_\de$ hull on $\Leb$.
Then, in Section \ref{s:all}, we prove that under \CH\ there is a monotone
Borel hull on $\Leb$, and prove partial results concerning $G_\de$ hulls. We
conclude the paper by collecting the open questions in Section \ref{s:conc}.

\section{Monotone hulls for nullsets}\lab{s:null}

%In this section we show that consistently there is a monotone $G_\de$ hull
%operation on $\iN$, hence consistently there is a monotone Borel hull
%operation on $\iN$, and also that consistently there is no monotone Borel hull
%operation on $\iN$, hence consistently there is no monotone $G_\de$ hull
%operation on $\iN$.

Recall that $\non(\iN) = \min\{|H| : H\su[0,1], \ H\notin\iN\}$, where $|H|$
denotes cardinality. In the sequel the cardinal $\ka$ is identified with its
initial ordinal, i.e.~with the smallest ordinal of cardinality $\ka$, and also
every ordinal is identified with the set of smaller ordinals. For the standard
set theory notation and techniques we use here see e.g.~\cite{Ku2} and
\cite{BJ}.

\bt
\lab{t:null-}
In a model obtained by adding $\om_2$ Cohen reals to a model
satisfying the Continuum Hypothesis (\CH) there is no monotone Borel hull
operation on $\iN$.
\et

\bp 
We need two well-known facts. Firstly, $\non(\iN) =\om_2$ in this model 
\cite{BJ}. Secondly, in this model there is no strictly increasing
(wrt.~inclusion) sequence of Borel sets of length $\om_2$ (this is proved in
\cite{Ku}, see also \cite{EK}). 

Assume that $\varphi: \iN\to\iB$ is a monotone hull operation. Choose
$H=\{x_\al:\al<\non(\iN)\} \notin \iN$, and consider
$\varphi(\{x_\beta:\beta<\al\})$ for $\al<\non(\iN)$. This is an increasing
$\om_2$ long sequence of Borel sets, which cannot stabilize, since then $H$
would be contained in a nullset. But then we can select a strictly increasing
subsequence of length $\om_2$, a contradiction.
\ep

The following is immediate.

\bcor
Under the same assumption there exists no monotone $G_\de$ hull
operation on $\iN$.
\ecor

\br
We will see in Remark \ref{r:o2} that the length $\om_2$ is optimal in the
sense that all shorter well-ordered chains have monotone $G_\de$ hulls.
\er

Recall that $\add(\iN) = \min\{|\iF| : \iF\su\iN, \ \bigcup\iF\notin\iN\}$ and
$\cof(\iN) = \min\{|\iF| : \iF\su\iN, \ \forall N\in\iN \ \exists F\in\iF
\textrm{ such that } N\su F\}$, and also that $\add(\iN) = \cof(\iN)$ is
consistent \cite{BJ} (note that e.g.~{\CH} implies this equality).

\bt
\lab{t:null+}
Assume $\add(\iN) = \cof(\iN)$. Then there exists a monotone $G_\de$ hull
operation on $\iN$.
\et

\bp 
Let $\{N_\al:\al<\cof(\iN)\}$ be a cofinal family in $\iN$, that is, $\forall
N\in\iN \ \exists \al<\cof(\iN) \textrm{ such that } N\su N_\al$. For every
$\al<\cof(\iN)$ define, using 
transfinite recursion, $A_\al = \,\textrm{a } G_\de \textrm{ hull of }
(\bigcup_{\beta<\al} A_\beta \cup N_\al)$. Clearly, $\{A_\al:\al<\cof(\iN)\}$ is a
cofinal increasing sequence of $G_\de$ sets. Now, for every $N\in\iN$ define
$\varphi(N) = A_{\al_N}$, where 
$\al_N$ is the minimal index for which $H\su A_{\al_N}$. It is easy to see that
$\varphi: \iN\to\iG_\de$ is a monotone hull operation.
\ep

The following is again immediate.

\bcor
Under the same assumption there exists a monotone Borel hull
operation on $\iN$.
\ecor

\section{Monotone hulls for all sets}\lab{s:all}

First we note (Statement \ref{s:equiv} below) that the title of this section is
justified, as there is no difference between working with measurable sets or
arbitrary sets. 

We need a well-known lemma first. 
Recall that the \emph{density topology} of $\RR$ consists of those measurable
sets that
have Lebesgue density $1$ at each of their points (see
e.g.~\cite{Ox}). Closure in this topology is denoted by $\overline{H}^d$, and
the term `hull' is used in the sense of Remark \ref{r:Polish}.

\bl
\lab{l:dens}
$\overline{H}^d$ is a hull of $H$ for every $H\su \RR$.
\el

\bp
Assume to the contrary that there exists 
a Lebesgue measurable set $L\su\RR$ with $\la(L)>0$ such that
$L\su \overline{H}^d\sm H$. Set $L_0 =
\{x\in L: x \textrm{ is a density point of } L\}$. By
the Lebesgue Density Theorem $L\sm L_0$ is a nullset, which easily implies that
$L_0\neq \eset$ is open in the
density topology. But $L_0\su\overline{H}^d$ is disjoint from $H$, a
contradiction.
\ep

\bst
\lab{s:equiv}
The existence of a monotone Borel/$G_\de$ hull operation on $\iP([0,1])$
%(in the sense defined in Remark \ref{r:rem} \ref{s})
is equivalent to the existence of a monotone
Borel/$G_\de$ hull operation on $\Leb$.
\est

\bp
On the one hand, the restriction to $\Leb$ of a monotone hull operation on
$\iP([0,1])$ is itself a monotone hull operation.

On the other hand, by the previous lemma there exists a monotone hull operation
$\psi:\iP([0,1])\to\Leb$ (note that $[0,1]$ is closed in the density topology).
Hence if $\varphi$ is a monotone hull operation on $\Leb$ then
$\varphi\circ\psi$ is a monotone hull operation on $\iP([0,1])$.
\ep

Theorem \ref{t:null-} immediately implies the following.

\bcor
In a model obtained by adding $\om_2$ Cohen reals to a model
satisfying \CH\  there is no monotone Borel or $G_\de$ hull
operation on $\Leb$.
\ecor

Now we turn to the positive results.

\bt
\lab{t:all+}
Assume \CH. Then there is a monotone Borel hull
operation on $\Leb$.
\et

Before we prove this theorem we need a few lemmas. In case $\iH=\iB$ the first
one is a special case of a well-known result about Borel liftings, but there
are no such results in case of $\iG_\de$. 

Let us denote by $A\Delta B$ the symmetric difference of $A$ and $B$.

\bl
\lab{l:lift}
(\CH) There exists a monotone map $\psi:\Leb\to\iG_\de$ such that
$\la(M\Delta\psi(M))=0$ for every $M\in\Leb$ and that $\la(M\Delta M')=0$
implies $\psi(M)=\psi(M')$ for every $M,M'\in\Leb$.
\el 

\bp
Let us say that $M,M'\in\Leb$ are \emph{equivalent}, if
$\la(M\Delta M')=0$. Denote by $[M]$ the equivalence class of $M$ and by
$\Leb/\iN$ the set of classes. We say that
$[M_1]\le [M_2]$ if there are $M_1'\in [M_1]$ and $M_2'\in [M_2]$ such that
$M_1'\su M_2'$. 

It is sufficient to define $\Psi: \Leb/\iN \to \iG_\de$ so that
$[M]\le [M']$ implies $\Psi([M])\su\Psi([M'])$ for every $M,M'\in\Leb$, and
that $\Psi([M])\in [M]$ for every $M\in\Leb$.

Enumerate $\Leb/\iN$ as $\{[M_\al]:\al<\om_1\}$. For every $\al<\om_1$ define
\[
\Psi([M_\al]) = \bigcap_{\substack{\beta<\al \\ [M_\beta] \ge [M_\al]}}
			\Psi([M_\beta]) \ \cap \,
            \Big( \textrm{a } G_\de \textrm{ hull of }
		    \bigcup_{\substack{\ga<\al \\ [M_\ga]\le [M_\al]}}
            \Psi([M_\ga]) \cup M_\al \Big).
\]
It is not hard to check that this is a $G_\de$ set such that $[M_\ga] \le
[M_\al] \le [M_\beta]$ implies $\Psi([M_\ga]) \su \Psi([M_\al]) \su
\Psi([M_\beta])$, and that $\Psi([M_\al])\in [M_\al]$, hence the construction
works.
\ep

\br %\mbox{} \par
\begin{enumerate}
\item
Actually we will not use the fact that $\psi$ is constant on the equivalence
classes.

\item
We do not know whether \CH\ is needed in this lemma, nor if \CH\ could be replaced by Martin's Axiom.
\end{enumerate}
\er

The following lemma is the only result we can prove for $\iB$ but not for
$\iG_\de$.

\bl\lab{l:vansub}
(\CH) There exists a monotone hull operation $\varphi:\iN\to\iB$ such that
\begin{enumerate}
\item
$\varphi(N \cup N') \su \varphi(N)\cup \varphi(N')$ for every $N,N'\in\iN$
(subadditivity),
\item\lab{2}
$\bigcup \{\varphi(N): N\su B, \, N\in\iN \} \sm B \in \iN$ for every $B\in\iB$.
\end{enumerate}
\el

\bp
Let $\{A_\al:\al<\om_1\}$ and $\al_N$ be as in Theorem \ref{t:null+} (note
that $\add(\iN)=\cof(\iN)=\om_1$ under \CH). Set $A_\al^* =
A_\al\sm\bigcup_{\beta<\al} A_\beta$. Enumerate $\iB$ as
$\{B_\al:\al<\om_1\}$ and for every $\al<\om_1$ define the countable set 
$$\iB_\al = \Big\{\bigcup_{i=0}^n B_{\beta_i}: n\in\NN,\, \beta_i<\al \ (0\le i\le n)\Big\}.$$
Note that every $\iB_\al$ is closed under finite unions.

Now define
\[
\varphi(N) = \bigcup_{\al\le\al_N} \Big( A_\al^* \cap
\bigcap_{\substack{B\in \iB_\al\\N\cap A_\al^*\su B}} B  \Big).
\]
This is clearly a disjoint union. It is easy to see that $\varphi$ is a
monotone Borel hull operation (note that $\varphi(N) \su A_{\al_N}$). 

For every $\al<\om_1$ define $\varphi_\al(N) = A_\al^* \cap \varphi(N)$
$(N\in\iN)$.
In order to check subadditivity, let $N,N'\in\iN$. We may assume
$\al_N\le\al_{N'}$, so clearly $\al_{N\cup N'} = \al_{N'}$.
It suffices to check that each $\varphi_\al$ is subadditive. 
If $\al>\al_N$ then actually $\varphi_\al(N \cup N') =
\varphi_\al(N')$, so we are done. Suppose now $\al\le\al_N$. Let $x\in A_\al^*$
such that $x\notin \varphi(N)\cup\varphi(N')$. Then there exist $B\supseteq
N\cap A_\al^*$ and $B'\supseteq N'\cap A_\al^*$ in $\iB_\al$ such that $x\notin B,
B'$. But then $B\cup B' \in \iB_\al$ witnesses that $x\notin \varphi(N\cup N')$
since $x\notin B\cup B'\supseteq (N\cup N')\cap A_\al^*$.

Finally, to prove \textit{\ref{2}} it is sufficient to show that $N\su B_\al$
implies 
$\varphi(N) \sm B_\al \su A_\al$ for every $N\in\iN$ and $\al<\om_1$. So let
$x\in \varphi_\beta(N)$ for some $\beta>\al$. We have to show $x\in
B_\al$. But this simply follows from the definition of $\varphi$ since
$B_\al\in\iB_\beta$. 
\ep

\bl
\lab{l:subeleg}
Let $\iH$ be either $\iB$ or $\iG_\de$.
Assume that there exists a monotone map $\psi:\Leb\to\iH$ such that
$\la(M\Delta\psi(M))=0$ for every $M\in\Leb$ and also that there exists a
monotone hull operation $\varphi:\iN\to\iH$ such that
\begin{enumerate}
\item
$\varphi(N \cup N') \su \varphi(N)\cup \varphi(N')$ for every $N,N'\in\iN$,
\item
$\bigcup \{\varphi(N): N\su H, N\in\iN \} \sm H \in \iN$
%$\cup_{N\su H} \varphi(N) \sm H \in \iN$
for every $H\in\iH$.
\end{enumerate}
Then $\varphi$ can be extended to a monotone hull operation
$\varphi^*:\Leb\to\iH$.
\el

\bp
We may assume that $\psi(N)=\emptyset$ for every $N\in\iN$ (by redefining
$\psi$ on $\iN$ to be constant $\emptyset$, if necessary).

Define
\[
\varphi^*(M) = \psi\big(M\big) \cup \varphi\big( M \sm \psi(M) \big) \cup
			\varphi\Big( \bigcup_{\substack{N\su\psi(M)\\
			\emptyset\neq N\in\iN}} \varphi(N) \sm \psi(M) \Big).
\]
%Note: "\emptyset\neq N" is only necessary in
%case \varphi(\emptyset)\neq\emptyset to ensure that \varphi^*=\varphi on
%negligible sets  
Clearly $\varphi^*(M)\in\iH$. As the union of first two terms contains $M$, we
obtain 
$M\su\varphi^*(M)$. Moreover, $\varphi^*(M)$ is a hull of $M$, since the first
term is equivalent to $M$ and the last two terms are nullsets. It is also
easy to see that $\varphi^*$ extends $\varphi$.

We still have to check monotonicity of $\varphi^*$. 
First we prove
\beq\lab{e}
N'\in\iN, \ M'\in\Leb, \ N'\su \psi(M') \Rightarrow \varphi(N') \su \varphi^*(M').
\eeq
Indeed, the case $N'=\eset$ is trivial to check, otherwise 
\begin{multline*}
\varphi(N') \su \bigcup_{\substack{N\su\psi(M')\\ \eset\neq N\in \iN}} \varphi(N) \su
\Big(
\bigcup_{\substack{N\su\psi(M')\\ \eset\neq N \in \iN}}
\varphi(N) \sm \psi(M')\Big) \cup \psi(M') \su \\
\su
\varphi\Big(\bigcup_{\substack{N\su\psi(M')\\ \eset\neq N \in \iN}}
\varphi(N) \sm \psi(M')\Big) \cup
\psi(M') \su \varphi^*(M'),
\end{multline*}
which proves \eqref{e}.
%holds if $N'\neq\eset$, otherwise \eqref{e} is trivial to check.

%First we prove this in the special case
%$N\su M$, $N\in\iN$, $M\in\Leb$. Using subadditivity we get 
%$\varphi^*(N)=\varphi(N) \su \varphi(N\cap \psi(M)) \cup \varphi(N\sm
%\psi(M))$. Now we estimate the two terms separately. Firstly, $\varphi(N\cap
%\psi(M)) \su \cup_{N'\su\psi(M)} \varphi(N') \su \left( \cup_{N'\su\psi(M)}
%\varphi(N') \sm \psi(M)\right) \cup \psi(M) \su \varphi\left(
%\cup_{N'\su\psi(M)} \varphi(N') \sm \psi(M) \right) \cup \psi(M) \su
%\varphi^*(M)$, where we also used that $\varphi$ is monotone. Secondly,
%$\varphi(N \sm \psi(M)) \su \varphi(M \sm \psi(M)) \su \varphi^*(M)$, so we
%are done with this special case.  

Let now $M\su M'$ be arbitrary elements of $\Leb$. We need to show that all
three terms of $\varphi^*(M)$ are contained in $\varphi^*(M')$.

Firstly, $\psi(M) \su \psi(M')$.

Secondly, define $N'=\big(M \sm \psi(M)\big) \cap \psi(M')$. Using
the subadditivity of $\varphi$ and then \eqref{e} we obtain
\begin{multline*}
\varphi\big(M \sm \psi(M) \big) \su
\varphi\Big(\big(M\sm \psi(M)\big) \cap \psi(M') \Big)
\cup \varphi\Big(\big(M \sm \psi(M)\big) \sm \psi(M') \Big) \su \\
\su \varphi\big(N'\big) \cup
\varphi\big(M' \sm \psi(M') \big) \su \varphi^*(M').
\end{multline*}

Thirdly, let
$$N'= \Big(\bigcup_{\substack{N\su\psi(M)\\ \eset\neq N\in \iN}}
\varphi(N) \sm \psi(M)\Big) \cap \psi(M').$$
Using the subadditivity of $\varphi$ and then \eqref{e} we obtain
\begin{multline*}
\varphi\Big( \bigcup_{\substack{N\su\psi(M)\\ \eset\neq N\in \iN}} \varphi(N)
\sm \psi(M) \Big)  
\su \\
%\su % overfull box miatt 
\varphi\Big( \Big(\bigcup_{\substack{N\su\psi(M)\\ \eset\neq N\in
\iN}}\varphi(N) \sm \psi(M)\Big) \cap \psi(M')\Big) 
%\left( [\cup_{N\su\psi(M)} \varphi(N) \sm \psi(M)] \cap \psi(M')\right) 
\cup
\varphi\Big( \Big(\bigcup_{\substack{N\su\psi(M)\\ \eset\neq N\in
    \iN}}\varphi(N) 
\sm \psi(M)\Big) \sm \psi(M') \Big)
%\left( [\cup_{N\su\psi(M)} \varphi(N) \sm \psi(M)] \sm \psi(M') \right)
\su \\
\su \varphi(N')
\cup
\varphi\Big( \bigcup_{\substack{N\su\psi(M')\\ \eset\neq N\in \iN}} \varphi(N)
\sm \psi(M') \Big)
\su
\varphi^*(M').
\end{multline*}
This concludes the proof.
\ep

Now we prove Theorem \ref{t:all+}.

\bp
Lemma \ref{l:lift} and Lemma \ref{l:vansub} show that in case of $\iH=\iB$ the
requirements of Lemma \ref{l:subeleg} can be satisfied, so the proof of
Theorem \ref{t:all+} is complete.
\ep

\br
\begin{enumerate}
\item
We remark that subadditive monotone maps are actually additive.

\item
The proof actually gives a monotone $F_{\si\de\si}$ hull. However, we do not
know whether a monotone $G_\de$ hull operation on $\Leb$ exists (Question
\ref{q}). Of course, in light of the previous theorem, under \CH, this is
equivalent to 
assigning $G_\de$ hulls only to the Borel (or $F_{\si\de\si}$) sets in a
monotone way. 

%Moreover, it is easy to see that it is also equivalent to
%assigning $G_\de$ hulls only to the $F_{\si\de}$ sets in a
%monotone way. ASSZEM MARHASAG!

\end{enumerate}
\er

\bq
Is there a monotone $G_\de$ hull operation on $\iB$? Or on $F_{\si\de\si}$? Or
on any other fixed Borel class e.g.~$\iF_\si$? (Of course $\iG_\de$ and the
simpler ones are not interesting.)
\eq

Our next goal is to prove Theorem \ref{t:chain}, the partial result we have
concerning monotone $G_\de$ hull operations on $\Leb$. 
It shows that it is not possible to prove in \ZFC\  the nonexistence of $G_\de$
hulls on $\iL$ along the lines of Theorem \ref{t:null-}, that is, only by 
considering long chains of sets.

%Cannot settle? Merthogy zfc-ben nem bizonyithato, hogy
%nincs omega_2-es G_delta novo sorozat? Vagy mire gondolsz?

\bt
\lab{t:chain}
%Assume $\add(\iN)=2^\om$.
Assume that there exists a monotone $G_\de$ hull operation $\psi$ on $\iN$
(which follows e.g.~from $\add(\iN)=\cof(\iN)$).
Let $\iC\su\iP([0,1])$ be a chain of sets, that is, for every $C,C'\in\iC$
either $C\su C'$ or $C'\su C$ holds. Then there exists a monotone $G_\de$ hull
operation on $\iC$.
\et

%From now on we will often use the obvious fact that if $M\su M'$ and
%$H$ is a hull of $M$ then $H\sm M'$ is a nullset.

\bp
By Lemma \ref{l:dens} we may assume that $\iC\su\iL$.

%We may also assume that $C\su [0,1]$ for every $C\in\iC$, since it is
%sufficient to construct the hulls separately in every $[n,n+1]$.
%%----ez nem kell, mert mar eleve [0,1]-ben vagyunk
Partition
$\iC$ into the intervals $\iI_r = \{ C\in\iC : \la(C)=r \}$. Let $R=\{r\in
[0,1] : \iI_r \neq \eset \}$, and fix an element $C_r \in \iI_r$ for every $r
\in R$. Well-order $R$ as $\{r_\al : \al<|R| \}$, and set $R_\al = \{r_\beta :
\beta<\al\}$. 

Now we define $\varphi(C_{r_\al})$ by transfinite recursion as follows. Fix
two countable sets $R_\al^- \su \{r\in R_\al : r < r_\al \}$ 
and $R_\al^+ \su \{r\in R_\al : r > r_\al \}$ 
so that $\forall r\in R_\al$,
$r<r_\al$ $\exists r' \in R_\al^-$ such that $r\le r'<r_\al$, and similarly,
$\forall r\in R_\al$, $r>r_\al$ $\exists r' \in R_\al^+$ such that
$r_\al<r'\le r$. (Note that $R_\al^-$ and $R_\al^+$ may be singletons or even
empty.) Set
\[
\varphi(C_{r_\al}) = \Bigg[ \textrm{a } G_\de \textrm{ hull of } \Bigg(
C_{r_\al} \cup \bigcup_{r\in R_\al^-} \varphi(C_{r}) \Bigg) \Bigg] \cap
\bigcap_{r\in R_\al^+} \varphi(C_{r}).
\]
It is easy to see that this is a monotone $G_\de$ hull operation on $\{C_r :
r\in R\}$.

We may assume that for the hull operation
$\psi$ we have $\psi(\eset)=\eset$.
Then we can define a monotone $G_\de$ hull
operation $\varphi_t$ on $\iI_t$ for each $t\in R$ as follows. Let
%By our assumption for every $t\in R$ we can define a monotone $G_\de$ hull
%operation $\varphi_t$ on $\iI_t$ as follows.
\[
\varphi_t(C) = \varphi(C_t) \cup \psi(C\sm C_t) \quad (C\in\iI_t).
\]  

For each $t\in R$ fix a countable set $R^{++}_t \su \{ r \in R : r > t \}$ 
so that $\forall
\,r\in R$, $r>t$ $\exists r'\in R^{++}_t$ such that $t<r'\le r$.
%But then we can finish the proof by setting
Set
\[
\varphi(C) = \varphi_{t}(C) \cap \bigcap_{r\in R^{++}_t} \varphi(C_r)
\]
for every $C\in\iI_t$ and every $t\in R$.
This is a proper definition since for $C=C_t$ this is just
an equality. It is easy to check that $\varphi(C)$ is a $G_\de$
hull of $C$ and that $\varphi$ is monotone.
\ep

Finally, we prove in \ZFC\ that rather long well-ordered chains have monotone
$G_\de$ hulls.

\bl
\lab{l:addig}
Let $\xi\le\add(\iN)$ and $\iC = \{ M_\al : \al<\xi \} \su \iP([0,1])$ be
such that $M_\al\su
M_\beta$ for every $\al\le\beta<\xi$. Then there exists a monotone $G_\de$
hull operation on $\iC$.
\el

\bp
By Lemma \ref{l:dens} we may assume that $\iC\su\iL$. 

By transfinite recursion define $A_\al$ to be a $G_\de$ hull of the set $M_\al
\cup \bigcup_{\beta<\al} (A_\beta \sm M_\al)$. Clearly every $A_\beta
\sm M_\al$ is a nullset, moreover there are $|\al|<\add(\iN)$ many of them,
hence $A_\al$ is a hull of $M_\al$, too.
\ep

Recall that $\ka^+$ is the successor cardinal of $\ka$ and also that every
$\xi<\ka^+$ has a cofinal (i.e.~unbounded) subset of order type at most
$\ka$.

\bt
\lab{t:add+ig}
Let $\eta<\add(\iN)^+$ and $\iC = \{ M_\al : \al<\eta \} \su \iP([0,1])$ be
such that
$M_\al\su M_\beta$ for every $\al\le\beta<\eta$. Then there exists a monotone
$G_\de$ hull operation on $\iC$.
\et

\bp
By Lemma \ref{l:dens} we may assume that $\iC\su\iL$.

We prove this lemma by induction on $\eta$.
Fix a cofinal subset $X \su \eta$ of order type $\xi\le\add(\iN)$ and also a
monotone $G_\de$ hull operation $\varphi_X:\{M_\al:\al\in X\}\to\iG_\de$ by
the previous lemma.  
Every complementary interval $[\beta,\ga)$ of $X$ (i.e.~every interval that is
maximal disjoint from $X$) is of order type $<\eta$, hence by the inductive
assumption there exists a monotone $G_\de$ hull operation
$\varphi_{[\beta,\ga)}:\{M_\al:\al\in [\beta,\ga)\}\to\iG_\de$. 
Also fix a measure zero $G_\de$ hull $H_{[\beta,\ga)}$ of $\bigcup_{\de<\beta,
\,\de\in X} \big(\varphi_X(M_\de)\sm M_\beta\big)$.
Now for every $[\beta,\ga)$ and every $\al\in [\beta,\ga)$ define
\[
\varphi(M_\al) = 
\left( H_{[\beta,\ga)} \cup \varphi_{[\beta,\ga)}(M_\al) \right) 
\cap 
\varphi_X(M_\ga),
\]
and also define $\varphi(M_\al) = \varphi_X(M_\al)$ for every $\al\in X$.
It is then easy to see that this is a monotone $G_\de$ hull operation on
$\iC$.
\ep

\br
\lab{r:o2}
As $\add(\iN)\ge\om_1$, we obtain that length $\om_2$ of the chain in the
proof of Theorem \ref{t:null-} was optimal. 
\er

\section{Concluding remarks and open problems}\lab{s:conc}

%As for $\subsetneqq$-preserving hulls, let us note that the case of $\Leb$ is
%easy (but we do not know the answer to the case of $\iN$).

%\bst
%There is no $\subsetneqq$-preserving monotone Borel hull on $\Leb$.
%\est

%\bp
%Let $C\su[0,1]$ be the Cantor set and let $B\su C$ be a Bernstein subset
%\cite{Ox}, that is, a set such that $B\cap F\neq\eset$ and $B\cap (C\sm
%F)\neq\eset$ for every uncountable closed set $F\su C$.  Then $C \sm A$ is
%countable for every Borel set $A$ containing $B$, as uncountable Borel sets
%contain uncountable closed sets \cite{Ke}.

%Clearly, $C \sm B$ is uncountable, so let $\{x_\al:\al<\om_1\}$ be distinct
%points of this set, then the strictly increasing chain $C_\al = ([0,1] \sm C)
%\cup B \cup \{x_\beta:\beta<\al\}$ cannot have a strictly monotone Borel hull
%$\varphi$, as already $\varphi(C_0)$ is of countable complement in $[0,1]$.
%\ep

First we show (in ZFC) that there are no \emph{strictly} monotone hulls of any
kind.

\bst
There is no $\subsetneqq$-preserving monotone Borel hull on $\iN$.
\est

\bp
It is well known that in every infinite set of size $\kappa$ there is a chain
(of subsets) of size greater than $\kappa$. Indeed, let $\la = \min\{\la' : 2^{\la'} >
\kappa\}$, 
and let us consider $X = \{ x \in 2^\la : \exists \alpha<\la \ \forall \beta \in
[\alpha, \la) \ x(\beta) = 0 \}$ ($2^\la$ is considered as the set of
  functions from $\la$ to $2 = \{0,1\}$).
Then $|X| = \kappa$, and it suffices to produce a $2^\lambda$-sized chain
of subsets of $X$. Let $<_{lex}$ denote the lexicographical ordering on
$2^\lambda$, and for $y \in 2^\lambda$ set $A_y = \{x \in X : x \le_{lex} y
\}$. Then $y <_{lex} y'$ implies $A_y \subsetneqq A_{y'}$, so
$\{A_y : y \in 2^\lambda\}$ is a chain of size $2^\lambda > \kappa$.

Hence the usual middle-third Cantor set (which is of measure zero) contains
a chain of size greater than continuum, but then the Borel hulls of the
elements of this chain form more than continuum many Borel sets, which is
impossible.
\ep

Now we pose a few somewhat vague problems, some of which may turn out to be
very easy.

\bq
It would be interesting to know what happens 
\begin{enumerate}
\item
if we look at the category analogue of Question \ref{q:main}, that is, when 
$\iN$ and $\Leb$ are replaced by the first Baire category (=meager) sets and
the sets with the property of Baire;

%\textbf{L\'etezik} monoton Borel burok a Baire-tulajdons\'ag\'u halmazokon! A
%biz ugyanez.

\item
if we require that our monotone hulls are translation or isometry invariant.

%\item
%if we replace $\su$ by $\subsetneqq$ in Question \ref{q:main}, that is, we
%require that strict containment is preserved.

\end{enumerate}
\eq

\bq
Does Theorem~\ref{t:all+} remain valid if we replace \CH\ by Martin's Axiom? That is, does there exist a monotone Borel hull operation on $\Leb$ if we assume Martin's Axiom?
\eq

%\bigskip

We now repeat the open questions of the paper for the sake of completeness.

\bq
Is there (in \ZFC) a monotone map $\psi:\Leb\to\iG_\de$ such that
$\la(M\Delta\psi(M))=0$ for every $M\in\Leb$? If yes, is there one such that
$\la(M\Delta M')=0$ implies $\psi(M)=\psi(M')$ for every $M,M'\in\Leb$?
\eq

\bq
Is there a monotone $G_\de$ hull operation on $\iB$? Or on $\iF_{\si\de\si}$? Or
on any other fixed Borel class e.g.~$\iF_\si$? (Of course $\iG_\de$ and the
simpler ones are not interesting.)
\eq

Let us conclude with the most important open question.

\bq
\lab{q}
Is it possible to assign $G_\de$ hulls to all (measurable) subsets of $[0,1]$ in a
monotone way?
\eq

\noindent
\textbf{Acknowledgement.} The authors are indebted to Mikl\'os Laczkovich and
Alain Louveau for some helpful comments. 
%We also gratefully acknowledge the support of \"Oveges Project of
%\includegraphics[width=1.3cm]{NKTH.eps} and
%\includegraphics[width=1.cm]{KPI.eps}. 

%~~~~~~~~~~~
%and 
%~~~~~~~~~. 
%\vspace{-0.83cm}
%\begin{figure}[h]~~~~~~~~~~~~~~~~~~~~~~~~~~~~~~~~~~~~~~~~~~~~~~~~~~~~~~~~~~
%~~~~~~~~~~~~~~~~~~~~~
%\end{figure} \vspace{-1.2cm}
%\begin{figure}[h]
% ~~~~~~~~~~~~~~~~~~~~~~~~~~~~~~~~~~~~~~~~~~~~~~~~~~~~~~~~~~~~~~~~~~~~~~~~
%~~~~~~~~~~~~~~~~~~~~~~~~
%\includegraphics[width=1.cm]{KPI.eps}
%\end{figure}

\bigskip

\noindent
\textsc{M\'arton Elekes} 

\noindent
\textsc{R\'enyi Alfr\'ed Institute of Mathematics} 

\noindent
\textsc{Hungarian Academy of Sciences} 

\noindent
\textsc{P.O. Box 127, H-1364 Budapest, Hungary}

\noindent
%\textit{Email:}
\verb+emarci@renyi.hu+

\noindent
%\textit{URL:}
\verb+http://www.renyi.hu/~emarci+

\bigskip

\noindent
\textsc{Andr\'as M\'ath\'e} 

\noindent
\textsc{E\"otv\"os Lor\'and University}

\noindent
\textsc{Department of Analysis}

\noindent
\textsc{P\'az\-m\'any P\'e\-ter s\'et\'any 1/c, H-1117 Budapest, Hungary}

\noindent
%\textit{Email:}
\verb+amathe@cs.elte.hu+

\noindent
%\textit{URL:}
\verb+http://amathe.web.elte.hu+

\end{document}